\documentclass{amsart}
\usepackage{amssymb}
\usepackage{graphicx}
\input xy
\xyoption{all}

\newenvironment{agrade}
  {\vspace*{0,2cm} \begin{center} { {\sc Acknowledgments }} \end{center}}%
  {\par \hfill }

\title{Convolution operators on Banach lattices with shift-invariant norms}
\author{Naz Miheisi}
\address{Department of Pure Mathematics, University of Leeds, Leeds LS2 9JT, U.K.}
\email{nazar@maths.leeds.ac.uk}
\date{October 2009}

\begin{document}

\maketitle

\begin{abstract}
Let $G$ be a locally compact abelian group and let $\mu$ be a complex valued regular Borel measure on $G$. In this paper we consider a generalisation of a class of Banach lattices introduced in \cite{joh}. We use Laplace transform methods to show that the norm of a convolution operator with symbol $\mu$ on such a space is bounded below by the $L^\infty$ norm of the Fourier-Stieltjes transform of $\mu$. We also show that for any Banach lattice of locally integrable functions on $G$ with a shift-invariant norm, the norm of a convolution operator with symbol $\mu$ is bounded above by the total variation of $\mu$.
\end{abstract}

\section{Introduction}
Let $G$ be a locally compact abelian (LCA) group. Let $M(G)$ denote the algebra of complex valued regular Borel measures on $G$. It is well known that for $L^p(G)$, a convolution operator with symbol $\mu\in M(G)$ (i.e an operator $f\mapsto f\ast\mu$ on $L^p(G)$) is bounded with norm less than or equal to $|\mu|(G)$. For $L^1(G)$ and $L^\infty (G)$, the norm of the convolution operator is equal to $|\mu|(G)$. Denote by $\widehat{G}$ the dual group of $G$ (the group of continuous homomorphisms $G\to\mathbb{T}$, where $\mathbb{T}$ denotes the circle group). It is also known that the norm of a convolution operator on $L^p(G)$ is greater than or equal to $|\hat{\mu}|_{\widehat{G}} = \sup_{\gamma\in\widehat{G}}|\hat{\mu}(\gamma)|$, where $\hat{\mu}$ $(=\mathcal{F}\mu)$ denotes the Fourier-Stieltjes transform of $\mu$. In this case, there is equality when the convolution operator acts on $L^2(G)$. For a comprehensive treatment of harmonic analysis on LCA groups see \cite{hew} or \cite{rud} for example.\\
In \cite{joh}, a class of function spaces defined on $\mathbb{R}_+$ - which we call Johansson spaces - was introduced and it was shown that the norms of convolution operators on these spaces has the same upper bound as for the $L^p$-spaces. We will describe these spaces in detail, and in the more general setting of LCA groups in section 2. Indeed, this upper bound holds for any function space defined on a discrete group provided that the norm is shift-invariant - that is, $\|f(x)\|=\|f(x-y)\|$ for fixed $y\in G$ - as is easily verified. Given this, it was conjectured that it holds for any space of locally integrable functions defined on $\mathbb{R}$ \cite{joh}, and then subsequently for any LCA group.\\
In section 2.1, we introduce Johansson spaces defined on LCA groups and give some basic results regarding the structure of these spaces. Then in section 2.2 we will restrict the domain of definition to being the positive cone (i.e all $x\ge0$) of an Archimedean ordered group. This allows us to use Laplace transform methods to obtain a lower bound for the norms of convolution operators on these spaces, which is the same as the lower bound for the $L^p$-spaces.\\
In section 3 we consider a much larger class of Banach lattices, which include both the $L^p$-spaces and Johansson spaces. We will show that provided the norm is shift-invariant, the norms of convolution operators on these spaces have the same upper bound as the $L^p$-spaces, thus answering the earlier conjecture positively for Banach lattices.

\subsection{Notation}
Throughout this paper $G$ will denote a $\sigma$-compact LCA group with identity $0$, Haar measure $m$ and dual group $\widehat{G}$. We write $\int f$ for $\int f\,\mathrm{d}m$. We denote the characteristic function of $E\subset G$ by $\chi_{_E}$. That is,
\[\chi_{_E}(x)=\left\{\begin{array}{lr} 
1 & \text{if} \quad x\in E\\
0 & \,\text{if} \quad x\notin E.
\end{array} \right.\]
$M(G)$ denotes the algebra of complex valued regular Borel measures on $G$ and $L_{loc}^1(G)$ denotes the Fr$\acute{\mathrm{e}}$chet space (cf. Lemma \ref{polish}) of complex valued measurable functions on $G$ such that \[\int_K|f|<\infty\] for every compact $K\subset G$. The topology on $L_{loc}^1(G)$ is that of $L^1$-convergence on compact subsets. For $f\in L_{loc}^1(G)$ and $\mu\in M(G)$ we define the convolution product $f\ast\mu$ by \[(f\ast\mu)(x)=\int_Gf(x-y)\,\mathrm{d}\mu(y).\] In addition, we define the shift operator $S_y:L_{loc}^1(G)\to L_{loc}^1(G)$ by $S_y f(x)=f(x-y)$.\\
Let $X$ be a normed linear space. By $X^*$ and $X_{[1]}$ we denote its topological dual and closed unit ball respectively. For $x\in X$ and $\phi\in X^*$, $\langle x,\phi\rangle$ denotes their dual pairing (i.e $\phi(x)$). We use the same notation for real or complex valued homomorphisms on $G$. If a complex valued function $f$ is bounded on some set $\Omega$, we define $|f|_\Omega:=\sup_{x\in\Omega}|f(x)|$. 

\section{Johansson spaces}
\subsection{Johansson spaces on LCA groups}
\newtheorem{d1}{Defintion}
\begin{d1} Let $K\subset G$ be a compact neighbourhood. For $1\leq p <\infty$, the Johansson space $J^p_K(G)$ consists of all functions $f\in L^1_{loc}(G)$ such that the norm \begin{displaymath} \|f\|_K=\sup_{x\in G}\left(\int_{K}|f(x-y)|^p\,\mathrm{d}m(y)\right)^{\frac{1}{p}} \end{displaymath} is finite. \end{d1}

\newtheorem{p1}[d1]{Proposition}
\begin{p1}\label{eqiv}
 Let $K,K'\subset G$ be compact neighbourhoods. Then $\|\cdot\|_K$ and $\|\cdot\|_{K'}$ are equivalent norms. \end{p1}

\begin{proof} Without loss of generality we can assume both $K$ and $K'$ to be neighbourhoods of the identity. Then $K\cap K'$ contains an open set $V$, so $\chi_{_V}\le\chi_{_K}$. Since $V$ contains the identity, the collection $\{V+a\}_{a\in K'}$ gives an open cover of $K'$. Since $K'$ is compact, there is a finite subcover $\{V+a_i\}_{i=1}^n$. Then we have \begin{displaymath} \chi_{_{K'}}\le\sum_{i=1}^n\chi_{_{V+a_i}}.\end{displaymath} Fix $x\in G$. Then it follows that
\begin{align} \int_G\chi_{_{K'}}(y)|f(x-y)|^p\,\mathrm{d}m(y) &\le\sum_{i=1}^n\int_G\chi_{_{V+a_i}}(y)|f(x-y)|^p\,\mathrm{d}m(y) \nonumber\\
&=\sum_{i=1}^n\int_G\chi_{_V}(y-a_i)|f(x-y)|^p\,\mathrm{d}m(y) \nonumber\\
&\le\sum_{i=1}^n\int_G\chi_{_K}(y-a_i)|f(x-y)|^p\,\mathrm{d}m(y) \nonumber\\
&=\sum_{i=1}^n\int_G\chi_{_K}(y)|f(x-y+a_i)|^p\,\mathrm{d}m(y) \nonumber\\
&\le\sum_{i=1}^n\|f\|_K^p. \nonumber 
\end{align}

And hence \begin{displaymath} \|f\|_{K'}\le n^{1/p}\|f\|_K. \end{displaymath}
Swapping $K$ and $K'$ in the preceding argument gives the result.
\end{proof}

Given the previous result, throughout the remainder of this paper we will avoid reference to the set $K$ and we will write $J^p(G)$ and $\|f\|$ for $J^p_K(G)$ and $\|f\|_K$ respectively. Proposition \ref{eqiv} also has the following interesting corollary.

\newtheorem{c2}[d1]{Corollary}
\begin{c2} (i) If $G$ is compact then $J^p(G)=L^p(G)$.\\
 (ii) If $G$ is discrete then $J^p(G)=\ell^\infty(G)$.
\end{c2}
\begin{proof} (i) Since all compact neighbourhoods $K$ give equivalent norms, we can take $K=G$ which gives $L^p(G)$.\\ (ii) If $G$ is discrete, singleton sets are open so we can take $K=\{0\}$ which gives $\ell^\infty(G)$. \end{proof}

For completeness we now present some basic results regarding the structure of $J^p(G)$. Proposition \ref{struct} answers a question that was asked in \cite{joh}.

\newtheorem{p3}[d1]{Proposition}
\begin{p3}\label{struct}
For $1\leq p <\infty$,\\ (i) $J^p(G)$ is a Banach space.\\
(ii) If $G$ is not compact, $J^p(G)$ is nonseparable.
\end{p3}

\begin{proof} Let the norm on $J^p(G)$ be induced by the compact neighbourhood $K$.\\
(i) Let $(f_n)_{n\in\mathbb{N}}\subset J^p(G)$ be a Cauchy sequence. There is a subsequence $(f_{n_k})_{k\in\mathbb{N}}$ such that $\|f_{n_{k+1}}-f_{n_k}\|\leq2^{-k}$. For each $m\in\mathbb{N}$, set
\begin{displaymath}
g_m=\sum_{k=1}^{m}|f_{n_{k+1}}-f_{n_k}|, \quad g=\sum_{k=1}^{\infty}|f_{n_{k+1}}-f_{n_k}|.
\end{displaymath}
For each $g_m$ we have
\begin{displaymath}
\|g_m\|\leq \sum_{k=1}^{m}\|f_{n_{k+1}}-f_{n_k}\|<1.
\end{displaymath}
Fix $y\in G$. By Fatou's lemma we see that
\begin{displaymath}
 \int_{y-K}\left(\liminf_{m\to\infty}g_m^p\right) \le \liminf_{m\to\infty} \int_{y-K}g_m^p \le \liminf_{m\to\infty}\|g_m\|^p \le 1,
\end{displaymath}
and hence $\|g\| \le 1$. Since $G$ is $\sigma$-compact it can be covered by countably many translates of $K$. So if $\|g\|<\infty$ then $g<\infty$ a.e, and so the series
\begin{displaymath}
 f(x)=f_{n_1}(x)+\sum_{k=1}^{\infty}\left(f_{n_{k+1}}(x)-f_{n_k}(x) \right)
\end{displaymath}
converges absolutely for almost all $x\in G$. Since
\begin{displaymath}
f_{n_1}(x)+\sum_{i=1}^{k-1}\left(f_{n_{i+1}}(x)-f_{n_i}(x)\right)\,=f_{n_k},
\end{displaymath}
we see that $f(x)=\lim_{k\to\infty}f_{n_k}(x)$ a.e. It now remains to show that $f$ is the norm limit of the sequence $(f_n)$.\\
Choose $\varepsilon>0$. There exists an $N\in\mathbb{N}$ such that $\|f_n-f_m\|<\varepsilon$ whenever $m,n\ge N$. For every $m>N$ and $x\in G$, Fatou's lemma shows that
\begin{displaymath}
 \int_{x-K}|f-f_m|^p \le \liminf_{k\to\infty}\int_{x-K}|f_{n_k}-f_m|^p \le \liminf_{k\to\infty}\|f_{n_k}-f_m\|^p \le \varepsilon^p.
\end{displaymath}
So we conclude $f\in J^p(G)$ with $f_n \to f$ in norm.\\

(ii) If $G$ is not compact, there is a sequence $(x_n)_{n\in\mathbb{N}}\subset G$ such that $\{x_n-K\}$ are pairwise disjoint. Let $\Omega\subset J^p(G)$ be the collection of all $f\in J^p(G)$ such that $f(G)\subset\{0,1\}$, $f$ is constant on each $x_n-K$ and $f(x)=0$ for almost all $x\in\cup_{n=1}^{\infty}(x_n-K)$. Clearly $\Omega$ is uncountable. For $f,g\in\Omega$ with $f\neq g$, there is an $N\in\mathbb{N}$ such that $f=1-g$ on $x_N-K$. From this it follows
\[\|f-g\|^p\ge\int_{x_N-K}|f-g|^p = \int_{x_N-K}1 = m(K).\]
Thus $J^p(G)$ cannot be separable.
\end{proof}

\subsection{Johansson spaces on ordered groups}
So far we have considered $J^P(G)$ for a general LCA group $G$. However, in this section we shall restrict $G$ to be an Archimedean ordered group (for a treatment of these see \cite{rud}). In fact, we will consider Johansson spaces defined only on the positive cone $G_+$ of $G$, where we define $G_+:=\{x\in G:x\ge0\}$. We define $G_-$ analogously. $J^p(G_+)$ is defined in precisely the same way as $J^p(G)$ but with $G_+$ replacing $G$ at each instance. We can consider this to be the closed subspace of $J^p(G)$ consisting of all $f\in J^p(G)$ with $f(x)=0$ for almost all $x\in G_-\setminus\{0\}$. The case $G=\mathbb{R}$ is particularly important due to its potential applications in linear control. We will not discuss these here, but the interested reader is referred to \cite{joh}.\\
In order to arrive at our main result about convolution operators on $J^p(G_+)$ we will require the machinery of Laplace transforms. We now give a brief discussion of this.\\
Let $\mathbb{C}^\bullet$ and $\mathbb{R}_+^\bullet$ denote the multiplicative groups of non-zero complex numbers and positive real numbers respectively. Define $\widetilde{G}:=$Hom$(G,\mathbb{C}^\bullet)$ - that is, $\widetilde{G}$ is the group of all continuous homomorphisms $G\to\mathbb{C}^\bullet$. Since $\mathbb{C}^\bullet=\mathbb{R}_+^\bullet\oplus\mathbb{T}$, we see that for each $x\in G$ and $\lambda\in\widetilde{G}$ we have $\langle x,\lambda\rangle=\langle x,\beta\rangle\langle x,\gamma\rangle$, where $\gamma\in\widehat{G}$ and $\beta\in$Hom$(G,\mathbb{R}_+^\bullet)$, and hence $\langle x,\lambda\rangle=\langle x,\gamma\rangle e^{\langle x,\alpha\rangle}$, where $\alpha\in$Hom$(G,\mathbb{R})$ (here we are considering the additive group of real numbers). Let $\Lambda\subset\widetilde{G}$ be the collection of all $\lambda\in\widetilde{G}$ with $|\langle x,\lambda\rangle|>1$ for all $x\notin G_-$. So for $\lambda\in\Lambda$ we have $\langle x,\lambda\rangle=\langle x,\gamma\rangle e^{\langle x,\alpha\rangle}$ with $\langle x,\alpha\rangle> 0$ for each $x\notin G_-$. Choose a compact neighbourhood $K\subset G$ and $y\in G$ such that the interval $[0,y]$ is contained in $K$. Let the norm on $J^p(G_+)$ be induced by this set. For $f\in J^1(G_+)$ and $\lambda\in\Lambda$,
\begin{align} \int_{G_+}|f(x)||\langle-x,\lambda\rangle|\,\mathrm{d}m(x) &=\int_{G_+}|f(x)|e^{-\langle x,\alpha\rangle}\,\mathrm{d}m(x) \\
&\le\sum_{n=0}^\infty\int_{ny-K}|f(x)|e^{-\langle x,\alpha\rangle}\,\mathrm{d}m(x) \nonumber\\
&\le\sum_{n=0}^\infty e^{-n\langle y,\alpha\rangle}\int_{ny-K}|f(x)|\,\mathrm{d}m(x) \nonumber\\
&\le\frac{\|f\|}{1-e^{-\langle y,\alpha\rangle}}. \nonumber \end{align}
Since $J^1(G_+)\supset J^2(G_+)\supset\dots$ \cite{joh}, we conclude that the integral (1) converges for all $f\in J^p(G_+)$ and $\lambda\in\Lambda$. This allows us to make the following definition.

\newtheorem{d3}[d1]{Definition}
\begin{d3} For $f\in J^p(G_+)$  (or $\mu\in M(G_+)$) we define the Laplace transform $\mathcal{L}:J^p(G_+)\to\mathbb{C}^\Lambda$ (resp. $M(G_+)\to L^\infty(\Lambda)$) by
\begin{displaymath}
\mathcal{L}f(\lambda)=\hat{f}(\lambda)=\int_{G_+}f(x)\langle-x,\lambda\rangle\, \mathrm{d}m(x) \qquad (\lambda\in\Lambda)
\end{displaymath} respectively,
\begin{displaymath}
\mathcal{L}\mu(\lambda)=\hat{\mu}(\lambda)=\int_{G_+}\langle-x,\lambda\rangle\, \mathrm{d}\mu(x) \qquad (\lambda\in\Lambda).
\end{displaymath}
\end{d3}

It is straightforward to verify that the Laplace transform of $\mu\in M(G_+)$ does indeed define a bounded function on $\Lambda$. For the case $G=\mathbb{R}$, $\Lambda$ corresponds to the right half plane $\mathbb{C}_+$ and this is just the ordinary Laplace transform as would be expected. The properties of the Laplace transform that we will be using are the following.

\newtheorem{p4}[d1]{Proposition}
\begin{p4}\label{prop} For $f\in J^p(G_+)$, $\mu\in M(G_+)$ and $\lambda\in\Lambda$,\\
(i) $\mathcal{L}(f\ast\mu)(\lambda)=\hat{f}(\lambda)\hat{\mu}(\lambda)$.\\
(ii) The linear functional $f\mapsto\hat{f}(\lambda)$, $J^p(G_+)\to\mathbb{C}$ is bounded for each $\lambda\in\Lambda$. \end{p4}

The proof of Proposition \ref{prop}(i) is identical to the one for the Fourier transform found in \cite{rud} and so we omit it here, and Proposition \ref{prop}(ii) follows immediately from the previous calculation showing the convergence of the integral (1). \qed \\

Since the functional $\mu\mapsto\hat{\mu}(\lambda)$ is a character on the Banach algebra $M(G_+)$ (under convolution) for every $\lambda\in\Lambda$, we can regard $\Lambda$ as a subset of $\Delta$, where $\Delta$ denotes the maximal ideal space of $M(G_+)$ (for a comprehensive treatment of Gelfand theory see \cite{dal}). The same is true for $\widehat{G}$. If we let $\mathcal{G}\mu$ denote the Gelfand transform of $\mu\in M(G_+)$, it follows that $\mathcal{L}\mu=\mathcal{G}\mu|_{\Lambda}$ and $\mathcal{F}\mu=\mathcal{G}\mu|_{\widehat{G}}$, so there is no ambiguity in using $\hat{\mu}$ to denote both $\mathcal{L}\mu$ and $\mathcal{F}\mu$. We will now show that $|\hat{\mu}|_{\widehat{G}}\le|\hat{\mu}|_{\Lambda}$.

\newtheorem{p5}[d1]{Proposition}
\begin{p5} For $\mu\in M(G)$, $|\hat{\mu}|_{\widehat{G}}\le|\hat{\mu}|_{\Lambda}$. \end{p5}

\begin{proof} First we will show that $\widehat{G}\subset\overline{\Lambda}$ when considered as subsets of $\Delta$ (which is equipped with the relative weak$^*$ topology).\\
Take $\gamma\in\widehat{G}$ and let $\alpha\in$Hom$(G,\mathbb{R})$ be a fixed homomorphism such that $\langle x,\alpha\rangle>0$ for all $x\notin G_-$. For each $n\in\mathbb{N}$ define $\lambda_n\in$Hom$(G,\mathbb{C}^\bullet)$ by
\[\langle x,\lambda_n\rangle = \langle x,\gamma\rangle e^{\frac{1}{n}\langle x,\alpha\rangle}.\]
It is apparent that $\lambda_n\in\Lambda$ for every $n\in\mathbb{N}$. We claim that $\lambda_n\to\gamma$ in $\Delta$. For each $\mu\in M(G_+)$ we have
\begin{align}
|\hat{\mu}(\lambda_n)-\hat{\mu}(\gamma)|&=\left\vert\int_{G}\langle -x,\gamma\rangle\left(1-e^{-\frac{1}{n}\langle x,\alpha\rangle}\right)\,\mathrm{d}\mu(x)\right\vert \nonumber\\
&\le\int_G\lvert 1-e^{-\frac{1}{n}\langle x,\alpha\rangle}\rvert\,\mathrm{d}|\mu|(x).\nonumber
\end{align}
Fix $\varepsilon>0$. Since $|\mu|$ is finite there exists $y\in G_+$ such that $|\mu|(\{x:x>y\})<\frac{\varepsilon}{2}$. Set
\[N=\left\lceil\frac{\langle y,\alpha\rangle}{\mathrm{log}\left(1-\frac{\varepsilon}{2|\mu|([0,y])}\right)}\right\rceil.\]
It is easily verified that for every $n\ge N$ we have
\[\lvert 1-e^{-\frac{1}{n}\langle x,\alpha\rangle}\rvert<\frac{\varepsilon}{2|\mu|([0,y])}\]
for each $x\in[0,y]$. Then for $n\ge N$
\begin{align}
|\hat{\mu}(\lambda_n)-\hat{\mu}(\gamma)|&<\int_{[0,y]}\lvert 1-e^{-\frac{1}{n}\langle x,\alpha\rangle}\rvert\,\mathrm{d}|\mu|(x) + \frac{\varepsilon}{2} \nonumber\\
&<\int_{[0,y]}\frac{\varepsilon}{2|\mu|([0,y])}\,\mathrm{d}|\mu|(x) + \frac{\varepsilon}{2} \nonumber\\
&=\varepsilon. \nonumber
\end{align}
Hence $\lambda_n\to\gamma$ in $\Delta$.\\
Since $\mathcal{G}\mu\in C_0(\Delta)$ (\cite{dal}, page 203 for example) we have $|\hat{\mu}|_{\Lambda}=|\hat{\mu}|_{\overline{\Lambda}}$, and given that $\widehat{G}\subset\overline{\Lambda}$, we conclude that $|\hat{\mu}|_{\widehat{G}}\le|\hat{\mu}|_{\Lambda}$.\end{proof}

In \cite{are} it was shown that if $\mu\ll m$, then $|\hat{\mu}|_{\Lambda}=|\hat{\mu}|_{\widehat{G}}$. This is a consequence of the natural correspondence between the maximal ideal space of $L^1(G)$ and $\widehat{G}$. However, if $G$ is not discrete, the maximal ideal space of $M(G)$ (or $M(G_+)$) is far more complicated and not easily characterised (see for example \cite{tay1}, \cite{tay2} or chapter 8 of \cite{gra}).\\

We can now give the main result of this section.

\newtheorem{t5}[d1]{Theorem}
\begin{t5}\label{thm} For $\mu\in M(G_+)$ define $T_\mu:J^p(G_+)\to J^p(G_+)$ by $T_\mu f=f\ast\mu$. Then $\|T_\mu\|\ge |\hat{\mu}|_{\widehat{G}}$. \end{t5}
\begin{proof} From Proposition \ref{prop}(ii) it follows that $e_\lambda:f\mapsto\hat{f}(\lambda)\in J^p(G_+)^\ast$ for each $\lambda\in\Lambda$. Given any $f\in J^p(G_+)$ and $\lambda\in\Lambda$ we have
\begin{displaymath} 
\langle f,T_\mu^*e_\lambda\rangle=\langle T_\mu f,e_\lambda\rangle=\hat{f}(\lambda)\hat{\mu}(\lambda)=\hat{\mu}(\lambda)\langle f,e_\lambda\rangle=\langle f,\hat{\mu}(\lambda)e_\lambda\rangle. \end{displaymath}
Hence each $\hat{\mu}(\lambda)$ is an eigenvalue of $T_\mu^*$ and consequently $\|T_\mu\|=\|T_\mu^*\|\ge|\hat{\mu}|_{\Lambda}\ge|\hat{\mu}|_{\widehat{G}}$. \end{proof}

Theorem \ref{thm} was already stated in \cite{joh} for the case $G=\mathbb{R}$. However, the proof contained an error: using the notation in \cite{joh}, it was stated that the functions $\bar{\gamma}_n\to 0$ as $n\to\infty$, which is false. As a result Theorem \ref{thm} corrects and extends the result. It is also worth noting that Theorem \ref{thm} is valid for any space of functions for which we can define a Laplace transform (or Fourier transform) and the functional $e_\lambda:f\mapsto\hat{f}(\lambda)$ is bounded for every $\lambda$.

\section{Banach Function Lattices}
\newtheorem{d8}[d1]{Definition}
\begin{d8} By a Banach function lattice (BFL) we mean a subspace $X$ of $L_{loc}^1(G)$ such that $(X,\|\cdot\|)$ is a Banach lattice (with the obvious lattice operations) satisfying the following condition:\\
 (1) $f\ast\nu\in X$ for each $\nu\in M(G)$ whenever $f\in X$.\\
We do not assume that the norm topology on $X$ is the same as the subspace topology inherited from $L_{loc}^1(G)$. \end{d8}
Note that Condition (1) above ensures that $S_y f\in X$ for each $y\in G$ whenever $f\in X$. This is because the shift operator $S_y$ can be viewed as a convolution operator with symbol $\delta_y$, where $\delta_y$ denotes the point mass at $y$.\\
In order to prove our main result about convolution operators on BFLs, we first need to establish some properties of the space $L_{loc}^1(G)$.

\newtheorem{l6}[d1]{Lemma}
\begin{l6}\label{polish} $L_{loc}^1(G)$ is a Polish space.\end{l6}
\begin{proof} Since $G$ is $\sigma$-compact, there is a collection $\{K_i\}_{i\in\mathbb{N}}$ of compact sets such that $G=\cup_iK_i$. Let $E_n=\cup_1^nK_i$. For each $n\in\mathbb{N}$ define the seminorm \begin{displaymath}\rho_n(f)=\int_{E_n}|f|. \end{displaymath}
To see that $L_{loc}^1(G)$ is completely metrisable, consider the metric $d$ given by
\begin{displaymath} d(f,g)=\sum_{n=1}^{\infty}\frac{2^{-n}\rho_n(f-g)}{1+\rho_n(f-g)}.\end{displaymath}
It is straightforward to verify that this is a complete metric and that it induces the correct topology.\\
Next we show separability. $L^1(E_n)$ is separable for each $n\in\mathbb{N}$ and so has a countable dense subset, $A_n$ say. Let $A=\cup_{n=1}^{\infty}A_n$. We claim that $A$ is dense in $L_{loc}^1(G)$.\\
By the separability of $L^1(E_n)$, for each $f\in L_{loc}^1(G)$ there exists a sequence $(f_{n,k})_{k\in\mathbb{N}}\subset A$ such that $\rho_n (f_{n,k}-f)<1/k$ for every $k\in\mathbb{N}$. Note that $\rho_n (f_{n,k}-f)<1/k$ implies $\rho (f_{m,k}-f)<1/k$ for all $m\le n$. Consider the sequence $(f_{k,k})_{k\in\mathbb{N}}$. For this sequence we have
\begin{align} d(f_{k,k},f)&=\sum_{n=1}^{\infty}\frac{2^{-n}\rho_n(f_{k,k}-f)}{1+\rho_n(f_{k,k}-f)} \nonumber\\
 &<\sum_{n=1}^{k}\frac{2^{-n}}{k+1}+\sum_{n=k+1}^{\infty}2^{-n} \nonumber\\
&< \frac{1}{k+1}+\sum_{n=k+1}^{\infty}2^{-n} \nonumber \end{align}
which clearly tends to $0$ as $k\to\infty$.
\end{proof}

\newtheorem{l7}[d1]{Lemma}
\begin{l7}\label{cts} The map $y\mapsto S_yf$, $G\to L_{loc}^1(G)$ is continuous for each $f\in L_{loc}^1(G)$. \end{l7}
\begin{proof} Here it is sufficient to show $S_{y_\alpha}f\to f$ in $L_{loc}^1(G)$ for every net $(y_\alpha)\subset G$ with $y_\alpha\to 0$.\\
Let $U\subset G$ be a compact, symmetric neighbourhood of the the identity and $K\subset G$ be an arbitrary compact set with $m(K)>0$. Fix $\varepsilon>0$ and $f\in L^1_{loc}(G)$. Since $f\in L^1_{loc}(G)$ and $K+U$ is compact, there exists a compactly supported, continuous function $g:G\to\mathbb{C}$ such that \[\int_{K+U}|f-g|<\frac{\varepsilon}{3}.\] So for $y\in U$ \[\int_{K}|S_y f-S_y g|<\frac{\varepsilon}{3}.\] Since $g$ is compactly supported, there is a symmetric neighbourhood $V$ of the identity such that $|g(x)-g(x-w)|<\varepsilon/3m(K)$ whenever $w\in V$. So if $y\in V\cap U$, \[\int_{K}|g-S_y g|<\int_{K}\frac{\varepsilon}{3m(K)}=\frac{\varepsilon}{3}.\] From this it follows that
\begin{align} \int_{K}|f-S_y f| &=\int_{K}|f-g+g-S_y g+S_y g-S_y f| \nonumber\\ &\le \int_{K}|f-g|+\int_{K}|g-S_y g|+\int_{K}|S_y g-S_y f|<\varepsilon. \nonumber \end{align}
\end{proof}

We now require some results about BFLs.

\newtheorem{l9}[d1]{Lemma}
\begin{l9}\label{bdd} Let $X$ be a BFL. Then for $K\subset G$ compact, the linear functional $X\to\mathbb{C}$, $f\mapsto\int_{K}f$ is bounded. \end{l9}
\begin{proof} For $K\subset G$ define $\mu:X\to\mathbb{C}$ by \begin{displaymath}\mu(f)=\int_{K}f .\end{displaymath}
Since $\mu$ is positive (i.e $\mu(f)\ge 0$ whenever $f\ge 0$) we have $f\ge g$ (equivalently $f-g\ge 0$) implies $\mu(f-g)\ge 0$ and hence $\mu(f)\ge\mu(g)$. We also have $\mu(|f|)\ge |\mu(f)|$.\\
Assume towards a contradiction that $\mu$ is not bounded. Then there exists a sequence $(g_n)_{n\in\mathbb{N}}\subset X$ such that $|\mu(g_n)|\ge1$ but $\|g_n\|<1/n^2$ for every $n\in\mathbb{N}$. Let $f_n=|g_n|$ so $f_n\ge0$, $\|f_n\|<1/n^2$ and $\mu(f_n)\ge|\mu(g_n)|\ge1$.\\
Let $f=\sum_{1}^{\infty}f_n$ which converges absolutely in $X$. For each $N\in\mathbb{N}$, \[ f\ge\sum_{1}^{N}f_n,\] and consequently \[\mu(f)\ge\sum_{1}^{N}\mu(f_n)\ge N, \]
which contradicts the fact that $f$ is locally integrable. \end{proof}

\newtheorem{c10}[d1]{Corollary}
\begin{c10}\label{inc} Every BFL is continuously included into $L_{loc}^1(G)$. \end{c10}
\begin{proof} For a BFL $X$, take a sequence $(f_n)_{n\in\mathbb{N}}\subset X$ such that $f_n\to0$. Then we have $|f_n|\to0$ by continuity of $|\cdot|$, and so by Lemma \ref{bdd} \[ \int_{K}|f_n|\to0\] for every compact $K\subset G$. Hence $f_n\to 0$ in $L_{loc}^1(G)$. \end{proof}

\newtheorem{l11}[d1]{Lemma}
\begin{l11}\label{borel} Let $X$ be a separable BFL. Then the map $y\mapsto S_y f$, $G\to X$ is Borel measurable for every $f\in X$. \end{l11}
\begin{proof} Let $i:X\hookrightarrow L_{loc}^1(G)$ be the natural inclusion, which we know to be continuous from Corollary \ref{inc}. Since $i$ is injective, we can define an inverse on its image $i^{-1}:i(X)\to X$, so for each $f\in X$ we have the following commutative diagram:\\
\[ \xymatrix{
X\ar@/^/[rr]^i & &i(X)\ar@/^/[ll]^{i^{-1}}\\
&G\ar[ul]^{y\mapsto S_y f}\ar[ur]_{y\mapsto S_y f} } \]
\\Since the composition of Borel maps is Borel and $y\mapsto S_y f$, $G\to i(X)\subset L_{loc}^1(G)$ is continuous by Lemma \ref{cts}, it only remains to show that $i^{-1}$ is a Borel map. To show this we first note that since $X$ is a separable Banach space it is trivially a Polish space, as is $L_{loc}^1(G)$ by Lemma \ref{polish}.\\
So for $E\subset X$ Borel, we have $(i^{-1})^{-1}(E)=i(E)$ which is the continuous injective image of a Borel set between Polish spaces and thus is itself Borel (\cite{kec}, Theorem 15.1). \end{proof}

We can now prove our main result.

\newtheorem{t12}[d1]{Theorem}
\begin{t12} Let $X$ be a BFL such that $\|S_y f\|=\|f\|$ for each $f\in X$ and $y\in\ G$. For $\mu\in M(G)$ define $T_\mu:X\to X$ by $T_\mu f=f\ast\mu$. Then $T_\mu$ is bounded and $\|T_\mu\|\le|\mu|(G)$. \end{t12}
\begin{proof} First we show that $T_\mu$ is bounded.\\
Define \begin{displaymath} \Gamma=\left\{f\mapsto\int_{K}f:K\subset G \quad\text{compact}\right\}\subset X^*. \end{displaymath}
For each $\psi\in\Gamma$ we have
\begin{align} |\langle T_\mu f, \psi\rangle| &=\left\vert\int_{K}\int_{G}f(x-y)\,\mathrm{d}\mu(y)\,\mathrm{d}m(x)\right\vert \nonumber\\
&\le\int_{G}|\langle S_y f,\psi\rangle|\,\mathrm{d}\mu \nonumber\\
&\le \|\psi\||\mu|(G)\|f\|. \nonumber \end{align}
So $\psi\circ T_\mu$ is bounded.\\
The second inequality follows from the boundedness of $S_y$ and $\psi$.\\
Take a sequence $(f_n)_{n\in\mathbb{N}}\subset X$ such that $f_n\to 0$ and $T_\mu f_n\to g$. By the boundedness of $\psi\circ T_\mu$ and $\psi$ it follows that $\langle T_\mu f_n, \psi\rangle\to0$ and $\langle T_\mu f_n, \psi\rangle\to\langle g,\psi\rangle$ respectively, and hence $\langle g,\psi\rangle=0$. Since $\Gamma$ separates the points of $X$, $g=0$. So by the closed graph theorem $T_\mu$ is bounded.\\
Next we show that for $X$ separable, $\|T_\mu\|\le|\mu|(G)$.\\
For $f\in X$ define $\lambda_f:G\to X$ by $\lambda_f(y)=S_y f$ which is Borel measurable by Lemma \ref{borel}. For each $\phi\in X^*$, the map $y\mapsto\langle\lambda_f(y),\phi\rangle$ is also Borel measurable since it is the composition of Borel maps, and given that $X$ is separable, Pettis's measurabilty theorem (\cite{die}, page 42) implies that $\lambda_f$ is strongly $\mu$-measurable. We also have
\begin{displaymath} \left\vert\int_{G}\|\lambda_f\|\,\mathrm{d}\mu\right\vert = \left\vert\int_{G}\|f\|\,\mathrm{d}\mu\right\vert = \|f\|\vert\mu(G)\vert <\infty, \end{displaymath}
so $\lambda_f$ is $\mu$-Bochner integrable. Define $R_\mu:X\to X$ by \begin{displaymath} R_\mu f=\int_{G}\lambda_f\,\mathrm{d}\mu. \end{displaymath}
For each $f\in X$ we have \begin{displaymath} \|R_\mu f\|=\left\|\int_{G}\lambda_f\,\mathrm{d}\mu\right\| \le\int_{G}\|\lambda_f\|\,\mathrm{d}|\mu| = |\mu|(G)\|f\|. \end{displaymath}
Hence $\|R_\mu\|\le|\mu|(G)$. We claim that $R_\mu=T_\mu$. For any $\psi\in\Gamma$ and $f\in X$ we have
\begin{align} \langle (R_\mu-T_\mu)f,\psi\rangle &=\left\langle\int_{G}\lambda_f\,\mathrm{d}\mu,\psi\right\rangle - \int_{K}\int_{G}f(x-y)\,\mathrm{d}\mu(y)\,\mathrm{d}m(x) \nonumber\\
&=\int_{G}\langle\lambda_f(y),\psi\rangle\,\mathrm{d}\mu(y) - \int_{K}\int_{G}f(x-y)\,\mathrm{d}\mu(y)\,\mathrm{d}m(x) \nonumber\\
&= \int_{G}\int_{K}f(x-y)\,\mathrm{d}m(x)\,\mathrm{d}\mu(y) - \int_{G}\int_{K}f(x-y)\,\mathrm{d}m(x)\,\mathrm{d}\mu(y) = 0. \nonumber \end{align}
The second equality followed directly from the fact that for $\phi\in X^*$ \begin{displaymath} \left\langle\int_{G}\lambda_f\,\mathrm{d}\mu,\phi\right\rangle = \int_{G}\langle\lambda_f,\phi\rangle\,\mathrm{d}\mu, \end{displaymath}
and the third equality followed from Fubini's theorem. Since $\Gamma$ separates the points of $X$ we see that $R_\mu=T_\mu$ and hence $\|T_\mu\|\le|\mu|(G)$.\\
Now we prove the general case. Let $X$ be nonseparable.\\
Choose $f\in X$ and define 
\begin{align} M_0 &=\mathrm{lin}\{f\} \nonumber\\ M_1 &=\mathrm{lin}(M_0\cup|M_0|\cup T_\mu M_0) \nonumber\\ M_2 &=\mathrm{lin}(M_1\cup|M_1|\cup T_\mu M_1) \nonumber\\ \vdots \nonumber \end{align}
where $|M_k|=\{|g|:g\in M_k\}$ and $T_\mu M_k=\{T_\mu g:g\in M_k\}$. Set $M=\overline{\cup_{k=1}^{\infty}M_k}$. For each $g\in M$ we have that $g=\lim_{n}g_n$ where $g_n\in M_{k_n}$ for some $k_n$. Consequently, $T_\mu g=\lim_{n}T_\mu g_n\in M$ and $|g|=\lim_{n}|g_n|\in M$ since $T_\mu$ and $|\cdot|$ are continuous and $M$ is closed. This shows that $M$ is a $T_\mu$-invariant Banach lattice.\\
We claim that $M$ is separable. Suppose that $M_k$ is separable for some $k\in\mathbb{N}$. Then $|M_k|$ and $T_\mu M_k$ are separable by the continuity of $|\cdot|$ and $T_\mu$, and so $\mathrm{lin}(M\cup|M_k|\cup T_\mu M_k)$ is separable. Since $M_0$ is separable it follows that $M_k$ is separable for every $k\in\mathbb{N}$ and hence $M$ is separable.\\
By considering the restriction of $T_\mu$ to $M$, the previous part of the proof can be applied. So for $T_\mu|_M^M:M\to M$ we have $\|T_\mu g\|\le|\mu|(G)\|g\|$ for each $g\in M$. In particular, $\|T_\mu f\|\le|\mu|(G)\|f\|$. Since $f$ was arbitrary, this holds for every $f\in X$.
\end{proof}

\begin{agrade}
The author would like to gratefully acknowledge the financial support of the UK Engineering and Physical Sciences Research Council (EPSRC), and the school of mathematics at the university of Leeds. He also wishes to thank Matthew Daws for many useful discussions and comments, as a result of which this work has been greatly improved, and Jonathan Partington for his ongoing guidance and support.
\end{agrade}

\end{document}